\documentclass{easychair}

\usepackage{rotating}
\usepackage{alltt}
\usepackage{pdflscape}

 \usepackage{algorithm2e}

\usepackage{amsmath, amsthm, amssymb}

\usepackage{dsfont}

\newcommand{\ff}{\mathds{F}}

\newtheorem{theorem}{Theorem}[section]   
\newtheorem{corollary}[theorem]{Corollary}     

\theoremstyle{definition}
\newtheorem{definition}[theorem]{Definition}   

\theoremstyle{remark}
\newtheorem{remark}[theorem]{Remark}        

\def\ds {\displaystyle}
\begin{document}

%
%
\title{Inferring Biologically Relevant Models: Nested Canalyzing Functions\thanks{The authors thank the referees for their valuable comments and suggestions.}}

%
\titlerunning{Inferring all nested canalyzing models}

\volumeinfo
	{h} 
	{1}                         
	{ANB10}    
	{1}                         
	{1}                         
	{1}                         

%
\author{
Franziska Hinkelmann\thanks{Hinkelmann was supported by a grant from the U.S. Army Research Office.}\\
Virginia Bioinformatics Institute\\
Virginia Tech\\
Blacksburg, VA 24060\\
\url{fhinkel@vbi.vt.edu}\\
\and
Abdul Salam Jarrah\thanks{The corresponding author.}\\
Department of Mathematics and Statistics\\
American University of Sharjah\\
Sharjah, United Arab Emirates\\
\url{ajarrah@aus.edu}\\
}

%
\authorrunning{Hinkelmann and Jarrah}

\maketitle

%
\begin{abstract}
Inferring dynamic biochemical networks is one of the main challenges
in systems biology. Given experimental data, the objective is to
identify the rules of interaction among the different entities of the network.
However, the number of possible models fitting the available data is huge and identifying a
biologically relevant model is of great interest. Nested canalyzing functions,
where variables in a given order dominate the function, have recently been proposed as a
framework for modeling gene regulatory networks.
Previously we described this class of functions as an algebraic toric variety.
In this paper, we present an algorithm that
identifies all nested canalyzing models that fit the given data. We demonstrate
our methods using a well-known Boolean model of the cell cycle in budding
yeast.
\end{abstract}

\section{Introduction}

Inferring dynamic biochemical networks is one of the main challenges
in systems biology. Many mathematical and statistical methods,
within different frameworks, have been
developed to address this problem, see \cite{reveng-review} for a review of some of these methods.
Starting from experimental data and known biological properties only,
the idea is to infer a ``most likely" model that could be used to generate
the experimental data. Here the model could have two parts. The first one is the
static network which is a directed graph showing the influence relationships among
the components of the network, where an edge from node $y$ to node $x$ implies that
changes in the concentration of $y$ could change the concentration of $x$.
The other part of the model
is the dynamics of the network, which describes how exactly
the concentration of $x$ is affected by that of $y$.
Due to the fact that biological networks are not well-understood and the available data
about the network is usually limited,
many models end up fitting the available information and the criteria
for choosing a particular model are usually not biologically motivated
but rather a consequence of the modeling framework.

A framework that has long been used for modeling
gene regulatory networks is \emph{time-discrete, finite-space
dynamical systems}. This includes Boolean networks \cite{Kauffman:69},
Logical models \cite{LM_thomas}, Petri nets \cite{Steggles}, and algebraic models \cite{LS}.
The latter is a straightforward generalization of Boolean networks
to multistate systems. Furthermore, in \cite{V-CJL}, it was shown that
logical models as well as Petri nets could be viewed and analyzed
as algebraic models. The inference methods we develop here are within the
algebraic models framework. To be self-contained,
 we briefly describe this framework
and state some of the known results that we need in this paper, see
\cite{LS,JLR,AB07,V-CJL} for more details.
Throughout this paper, we will be talking about gene
regulatory networks, however, the methods apply for biochemical networks in general.

Suppose that the gene regulatory network that we want to infer has $n$
genes and that we have a set $D$ of $r$ state transition pairs
$(\mathbf s_j, \mathbf t_j)$, $j=1,...,r$. The
input $\mathbf s_j$ and the output $\mathbf t_j$ are $n$-tuples
of 0 and 1 encoding the
state of genes $x_1,\dots,x_n$. Real time data points are not Boolean
but could be discretized (and in particular, could be made Boolean)
using different methods \cite{DVML}.
The goal is to find a model
\[
f = (f_1, f_2, \ldots ,f_n): \ff_2^n\longrightarrow \ff_2^n
\]
such that, for $j=1,\dots, r$,
\[
f(\mathbf s_j) = (f_1(\mathbf s_j),\ldots , f_n(\mathbf s_j)) = \mathbf t_j.
\]
Notice that, since any function over a finite field is a polynomial,
each $f_i$ is a polynomial.
An algorithm that finds all models $f$ is presented in \cite{LS}.
This is done by identifying, for each gene $i$,
the set of all possible functions for $f_i$.
This set can be represented as
the coset $f + I$, where $f$ is a particular such function and
$I\subset \ff_2[x_1, \ldots ,x_n]$ is the ideal of all Boolean polynomials
that vanish on the input data set, that is,
$I = \mathds{I}(\{\mathbf s_1,\dots, \mathbf s_r\})$.
The algorithm in \cite{LS} then proceeds to find a particular model
from the model space $f+I$. The chosen
model, which is the normal form of $f$ in the ideal $I$,
depends on the term ordering used in the Gr\"obner bases
computation. So different ordering of the variables (genes) might lead
to the selection of different models.
This presents a problem as the term ordering, which is a needed for
computational reasons, clearly influence the
model selection process.

Several modifications have since been
presented to address this problem. For example, in \cite{DJLS},
using the Gr\"obner fan of the ideal $I$, the authors developed
a method that produces a probabilistic model using all possible normal
forms. Other improvements on this algorithm can be found in \cite{JLSS,SJLS}.

Another approach toward improving the model selection process is
by restricting the model space $f+I$ by
requiring not only that the chosen model fits the data but also satisfies
some other conditions, such as
its network being sparse or scale-free, the polynomials $f_i$ being monomials,
the dynamics of the model having some desirable properties such as
fixed points are the only limit cycles (that is, starting from any initialization,
the model always reaches a steady state), or that the model is robust
and stable which could roughly mean that the number of attractors in the phase space is small.
In a nutshell, some but not all functions in the model space $f+I$ are biologically relevant
and hence restricting the space to only relevant models will improve the model selection process.

By desiring a particular property, several classes of functions have been proposed
as biologically relevant functions such as biologically meaningful rules \cite{Raey_2002},
certain post classes of Boolean functions
have been studied in \cite{shmul_2003}, and chain functions in
\cite{chain}, to name few.
Another class of Boolean functions, which was introduced by
S. Kauffman \emph{et al.} \cite{Kauffman:04}, is called (nested) canalyzing functions (NCF), where an input to a single variable exclusively
determines the value of the function regardless of the values of all other variables.
This is a natural characterization of ``canalisation" which was introduced by geneticist C. H. Waddington
\cite{wad} to represent the ability of a genotype to produce the
same phenotype regardless of environmental variability.
Indeed, known biological functions have been shown to be canalyzing
\cite{harris,wilhelm_2007}, and Boolean nested canalyzing networks to be robust and stable \cite{Kauffman:03,Kauffman:04,wilhelm_2007}.

For the purpose
of restricting the model space $f+I$ of all Boolean polynomial models to NCFs only,
we previously studied nested canalyzing functions, gave necessary and sufficient
conditions on the coefficients of a boolean polynomial function to be nested canalyzing, and showed that NCFs are
nothing but unate cascade functions \cite{JLR}.
Furthermore, in \cite{AB07}, the class of
all nested canalyzing functions is parameterized as the rational points of an
affine algebraic variety over the algebraic closure of $\ff_2$.
This variety was shown to be toric, that is,
defined by a collection of binomial polynomial
equations. In this paper we present an algorithm that restricts the model
space to only nested canalyzing functions by identifying
all NCFs from the model space $f+I$ that fit the given data set.

In the next section we briefly recall some definitions and results from
\cite{JLR,AB07}. Our algorithm is presented in Section \ref{sec:algorithm},
and its implementation in Singular is discussed in Section \ref{sec:sing-lib}.
Before we conclude this paper, we demonstrate the algorithm in Section
\ref{cell-cycle-sec},
where we identify all nested canalyzing models for the cell cycle in budding yeast using
time course data from the Boolean model in \cite{Tang:2004}.

\section{Nested Canalyzing Functions: Background}

We recall some of the definitions and the results from \cite{JLR,AB07} that we need
to make this presentation self-contained. Throughout this paper,
when we refer to a function on $n$ variables, we mean that $h$
depends on all $n$ variables, that is, for $i=1,\dots, n$,
there exists $(a_1,\dots,a_n) \in \ff_2^n$ such that
$h(a_1,\dots,a_{i-1},a_i,a_{i+1},\dots,a_n) \neq
h(a_1,\dots,a_{i-1},1+a_i,a_{i+1},\dots,a_n)$.
\begin{definition}
Let $h$ be a Boolean function on $n$ variables, i.e., $h:\mathbb \ff_2^n \rightarrow \mathbb \ff_2$.
\begin{itemize}
\item The function $h$ is a \emph{nested canalyzing function} (NCF) with respect to a permutation $\sigma$ on
the $n$ variables, canalyzing input value $a_i$ and canalyzed output value
$b_i$, for $i=1,\dots, n$, if it can be represented in the form
\begin{equation} \label{ncf kauff}
h(x_1,x_2,\ldots, x_n) =
\begin{cases}
    b_1 & ~{\rm if}~ x_{\sigma(1)} = a_1, \\
    b_2 & ~{\rm if}~ x_{\sigma(1)} \ne a_1 ~{\rm and}~ x_{\sigma(2)} = a_2, \\
    b_3 & ~{\rm if}~ x_{\sigma(1)} \ne a_1 ~{\rm and}~ x_{\sigma(2)} \ne a_2 ~{\rm and}~ x_{\sigma(3)} = a_3, \\
    \vdots & \hspace{1cm} \vdots \\
    b_n & ~{\rm if}~ x_{\sigma(1)} \ne a_1 ~{\rm and}~ \cdots ~{\rm and}~ x_{\sigma(n-1)} \ne a_{n-1} ~{\rm and}~ x_{\sigma(n)} = a_n, \\
    \overline{b_{n}} & ~{\rm if}~ x_{\sigma(1)} \ne a_1 ~{\rm and}~ \cdots ~{\rm and}~ x_{\sigma(n)} \ne a_n.
\end{cases}
\end{equation}
\item The function $h$ is \emph{nested canalyzing} if $h$ is nested canalyzing
with respect to some permutation $\sigma$, canalyzing input values $a_1,\dots, a_n$ and canalyzed output values
$b_1,\dots,b_n$, respectively.
\end{itemize}
\end{definition}

\begin{remark}
The definition above has been generalized to multistate functions in \cite{david},
where it is also shown that the dynamics of these functions are similar to their Boolean counterparts.
In \cite{ode-unate}, the authors introduce what they called \emph{kinetic models with unate structure},
which are continuous models having the canalization property, and they presented an algorithm
for identifying such models.
\end{remark}

Using the polynomial form of any Boolean function,
the ring of Boolean functions
is isomorphic to the quotient ring $R = \ff_2[x_1,\dots,x_n]/J$,
where $J = \langle x_i^2-x_i : 1 \leq i \leq n \rangle$. Indexing
monomials by the subsets of $[n]:=\{1,\ldots ,n\}$ corresponding
to the variables appearing in the monomial, the
elements of $R$ can be written as
\begin{equation*}
R = \{\ds \sum_{S \subseteq [n]} c_S \prod_{i \in S} x_i \, : \, c_S
\in \ff_2\}.
\end{equation*}
As a vector space over $\ff_2$, $R$ is isomorphic to $\ff_2^{2^n}$
via the correspondence
\begin{equation}\label{corresp}
R \ni \ds \sum_{S \subseteq [n]} c_S \prod_{i \in S} x_i
\longleftrightarrow (c_{\emptyset},\dots, c_{[n]}) \in
\ff_2^{2^n}.
\end{equation}

The main result in \cite{JLR} is the identification of the set of nested
canalyzing functions in $R$ with a subset $V^{ncf}$ of
$\ff_2^{2^n}$ by imposing relations on the coordinates of its
elements.

\begin{definition}
Let $\sigma$ be a permutation of the elements of the set $[n]$. We
define a new order relation $<_{\sigma}$ on the elements of $[n]$
as follows: $ \sigma (i) <_{\sigma} \sigma (j)$ if and only if $i
< j$. Let $r_S^{\sigma}$ be the maximum element of a nonempty
subset $S$ of $[n]$ with respect to the order relation
$<_{\sigma}$. For any nonempty subset $S$ of $[n]$, the {\it
completion of S with respect to the permutation $\sigma$}, denoted
by $[r_S^{\sigma}]$, is the set $[r_S^{\sigma}] = \{ \sigma(1),
\sigma(2), \ldots, \sigma(r_S ^\sigma) \}$.

Note that, if $\sigma$ is the identity permutation, then the
completion is $[r_S]$ := $\{1,2,\ldots,r_S\}$, where $r_S$ is the
largest element of $S$.
\end{definition}

\begin{theorem} \label{coeff rel sigma}
Let $h \in R$ and let $\sigma$ be a permutation of the set $[n]$.
The polynomial $h$ is nested canalyzing with respect to $\sigma$,
input value $a_i$ and corresponding output value
$b_i$, for $i=1,\dots, n$, if and only if $c_{[n]} = 1 $
and, for any proper subset $S \subseteq [n]$,
\begin{equation} \label{coeff formula sigma}
c_S = c_{[r_S^{\sigma}]} \prod_{\sigma (i) \in [r_S^{\sigma}] \backslash S} c_{[n] \backslash \{\sigma (i)\}}.
\end{equation}
\end{theorem}

\begin{corollary}\label{variety-sigma}
The set of points in $\ff_2^{2^n}$ corresponding to the set of all nested canalyzing functions
with respect to a permutation $\sigma$ on $[n]$, denoted by $V_{\sigma}^{ncf}$, is defined by
\begin{equation}
V_{\sigma}^{ncf} = \{(c_\emptyset,\dots,c_{[n]}) \in \ff_2^{2^n} :
c_{[n]}=1, \, c_S = c_{[r_S^{\sigma}]} \prod_{\sigma (i) \in
[r_S^{\sigma}] \backslash S} c_{[n] \backslash \{\sigma (i)\}} \mbox
{, for } S \subseteq [n]\}.
\end{equation}
\end{corollary}

It was shown in \cite{AB07} that $V_{\sigma}^{ncf}$ is an algebraic variety,
and its ideal
$\mathds{I}(V_\sigma^{ncf})$ is a binomial prime ideal in the polynomial ring
$\overline{\ff}_2[\{c_S : S \subseteq [n]\}]$, where
$\overline{\ff}_2$ is the algebraic closure of $\ff_2$. Namely,
\[
I_\sigma = \mathds{I}(V_\sigma^{ncf}) = \langle c_{[n]} -1, \, c_S - c_{[r_S^\sigma]}
\prod_{\sigma(i) \in [r_S^\sigma]\setminus S} c_{[n]\setminus
\{\sigma(i)\}} : S \subset [n]\rangle.
\]

Furthermore, the variety of all nested canalyzing functions is
\begin{eqnarray*}
V^{ncf} &=& \bigcup_\sigma V_\sigma^{ncf}
\end{eqnarray*}
and its ideal is
\[
\mathds{I}(V^{ncf}) = \bigcap_{\sigma} I_\sigma.
\]

In the next section, we identify the set $f+I$ with the rational points in
an algebraic affine variety. This will allow us to identify all nested canalyzing
functions in the model space $f+I$.

\section{Nested Canalyzing Models} \label{sec:algorithm}
Recall that we are given the data set
$D=\{(\mathbf s_1, \mathbf t_1), \dots,(\mathbf s_r, \mathbf t_r) \} \subset \ff_2^n \times \ff_2^n$.
The model space could be presented by the set $f+I$ where, $f=(f_1,\dots,f_n)$ and, for $i=1,\dots,n$,
\begin{equation}\label{equ:f}
f_i(x_1,\dots,x_n)=\sum_{j=1}^{r} t_{j,i} \prod_{e=1}^n (1-(x_e - s_{j,e})).
\end{equation}
In particular,
$f_i$ is a polynomial that interpolates the data for gene $i$ 
and $I$ is the \emph{ideal of points} of $\{\mathbf s_1,\dots,\mathbf s_r\}$.
Furthermore, the ideal $I$ is a principal ideal in the ring $R/J$:
\begin{eqnarray}
I &=& \mathds{I}(\{\mathbf s_1,\dots, \mathbf s_r\}) \\
  &=& \bigcap_{j=1}^r \mathds{I}(\{\mathbf s_j\})\\
  &=& \bigcap_{j=1}^r \langle x_1 - s_{j,1}, \ldots, x_n- s_{j,n} \rangle \\
  &=& \bigcap_{j=1}^r \langle 1 - \prod_{e=1}^n (1-(x_e-s_{j,e}))  \rangle \\
  &=& \langle \prod_{j=1}^r( 1-\prod_{e=1}^n (1-(x_e-s_{j,e}))) \rangle.  \label{equ:I}
\end{eqnarray}
Now a polynomial
$h \in f_i +I$ if and only if $h = f_i + g(x_1,\dots,x_n) \prod_{j=1}^r( 1-\prod_{e=1}^n (1-(x_e-s_{j,e})))$,
for some polynomial $g$, say $g = \sum_{H\subseteq [n]} b_H \prod_{i \in H} x_i$.
By expanding the right-hand side and collecting terms, we get that
$h  = \sum_{S \subseteq [n]} W_S(b_H,\mathbf s_j, \mathbf t_j) \prod_{l \in S} x_l$, where, for $S \subseteq [n]$,
the coefficient $W_S(b_H,\mathbf s_j, \mathbf t_j)$ is determined by
$b_H,\mathbf s_j, \mathbf t_j$ for all $H \subseteq [n]$ and $j=1,\dots,r$.

The proof of the following theorem follows directly from Theorem 2.4.2 in \cite{AdamsLost}.
\begin{theorem}
Consider the ring homomorphism
\[
\Phi : \overline{\ff}_2[\{c_S : S\subseteq [n]\}] \longrightarrow \overline{\ff}_2[\{b_H : H \subseteq [n]\}]
\]
given by, for $S \subseteq [n]$,
\begin{eqnarray*}
    c_{S} &\mapsto& W_S(b_H,\mathbf s_j, \mathbf t_j).
\end{eqnarray*}
Then $\ker(\Phi)$ is the ideal of all polynomials that fit the data set $D$. In particular,
the rational points in the variety $\mathds{V}(\ker(\Phi))$ is the set of all
models that fit
the data set $D$, namely $f+I$.
\end{theorem}
Since the ideal of all NCFs is $\mathds{I}(V^{ncf})$, the following corollary is straightforward.
\begin{corollary}
The ideal of all nested canalyzing functions that fit the data set $D$ is $\mathds{I}(V^{ncf})+\ker(\Phi)$.
\end{corollary}

\begin{remark}
It is clear that the model space of Boolean functions is huge, since the
number of monomials grows exponentially in the number of variables.
For example, if a function has 5 inputs, there are $2^5 =32$ different monomials in $5$
variables, and hence $2^{32}=4,294,967,296$ different Boolean functions.
This clearly shows that a search for NCFs inside the model space is
 computationally not feasible, which justifies the need for algorithms like the one above.
\end{remark}
\section{Algorithm}\label{sec:sing-lib}
In this section we present an algorithm for identifying all nested canalyzing models from the model space of a given data set.
\subsubsection*{Input}
A wiring diagram, i.e., a square matrix of dimension $n$, describing the
influence relationships among the $n$ genes in the network.
For each variable $x_i$, a table consisting of the rows
$(s_{j, i_1}, \ldots, s_{j, i_s}, t_{j,i})$, $j = \{ 1, \ldots, r\}$, where
$i_1, \ldots, i_s$ are the indices of the genes that affect $x_i$, as specified in the wiring diagram.
\subsubsection*{Output}
For each variable, the complete list of all nested canalyzing functions interpolating the
given data set on the given wiring diagram. A function is in the output if it is nested canalyzing in at least one variable order.
If needed, the code can easily be modified to find only nested canalyzing functions
of a particular variable order.
\subsubsection*{Algorithm}
It is a well known fact, that a Gr\"obner basis for the kernel of $\Phi$ is
a basis for $\langle c_S - W_S: S\subseteq [n]\} \rangle$ intersected with the ring
$\overline{\ff}_2[\{c_S : S\subseteq [n]\}]$ \cite[Theorem 2.4.2]{AdamsLost}
Using a similar notation as above, the algorithm is outlined as follows:
\begin{tabbing}
Use ring $\mathbb F_2[x_1, \ldots, x_n, b_S, c_S: S \subseteq [n]\} ]$\\
Define $\mathds{I}(V^{ncf})$ as ideal in $\mathbb F_2[ c_S : S \subseteq [n]\}]$\\
Define $h = \sum_{H\subseteq [n]} b_H \prod_{i \in H} x_i$\\
Define $q = \sum_{H\subseteq [n]} c_H \prod_{i \in H} x_i$\\
Compute the polynomial $p$ that generates $I$ as in (\ref{equ:I}) \\
For \= each variable $x_i$ do\\
\>  1. Compute $f_i$ as in (\ref{equ:f})\\
\>  2. Let $g = f_i + h*p$; its coefficients are the same as $W_S$ above\\
\>  3. Compute a Gr\"obner basis $G$ for the ideal generated by the coefficients of $g-q$ using any \\
\>     elimination order to eliminate all $b_S$ from $G$\\
\>  4. Concatenate generators of $G$ and $\mathds{I}(V^{ncf})$\\
\>  5. Compute the primary decomposition of $G + \mathds{I}(V^{ncf})$ to obtain necessary and sufficient conditions \\
\>     on the coefficients of all NCFs fitting the data set $D$\\
End
\end{tabbing}
An implementation of this algorithm is available as a Singular library
\cite{Hinkelmann, singular}.
\section{Application: Inferring the Cell Cycle Network in Budding Yeast}\label{cell-cycle-sec}
Li \emph{et al.} \cite{Tang:2004} constructed a Boolean threshold model of the cell cycle in budding yeast.
The network of the model has the key known regulators of the cell cycle process and
the known interactions among these regulators in
the literature. The Boolean function at each node, however, is a threshold function,
which is completely determined by the numbers of active activators and active inhibitors,
and is not necessarily biologically motivated.
However, this model captures the known features of the global dynamics of the cell cycle,
it is robust and stable, and the trajectory of the known cell-cycle sequence
is a stable and attracting trajectory as it has 1764 states out of the total number of 2048 states.
The remaining states are distributed into 6 small trajectories.

In this Section, we use the time course corresponding to the biological cell-cycle sequence, see Table \ref{table:timecourse},
to infer nested canalyzing models of the cell cycle. That is, assuming the same wiring diagram
as the threshold model above, we use our algorithm to identify, for each gene in
the network, all nested canalyzing functions that fit the
cell cycle sequence.

We start by describing the network. It
consists of 11 proteins and a start signal. The proteins are members of the
following three classes:
cyclins (Cln1,2, Cln3, Clb1,2, and Clb5,6), inhibitors, degraders, and
competitors of cyclin complexes (Sic1, Cdh1, Cdc20 and Cdc14), and transcription
factors (SBF, MBF, Swi5, and Mcm1/SFF).
\begin{figure}[htb]
  \includegraphics[width=.9\textwidth]{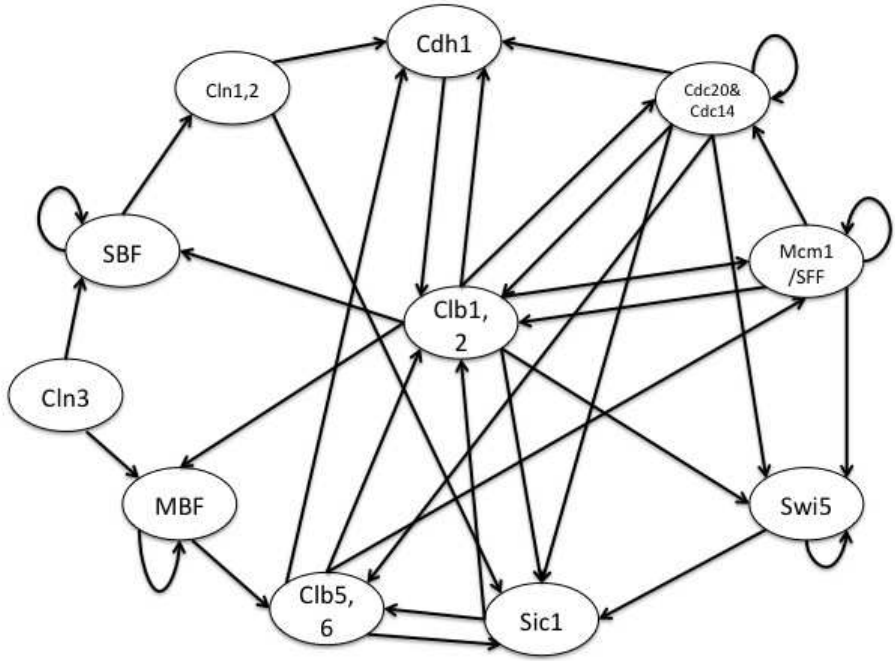}\\
  \caption{The simplified cell cycle network in budding yeast, which is based on the model in \cite{Tang:2004}.}
  \label{fig:dep}
\end{figure}
This simplified network
(Figure \ref{fig:dep}) is almost identical to the network in \cite{Tang:2004} where the only
difference is that we do not force self-degradation, as it was added to some nodes
in the network because they did not have inhibitors, but without biological justification \cite{Tang:2004}.
Furthermore, we do not impose
activation or inhibition in the network. As we do not use threshold functions
but more general boolean functions, a variable can both increase and decrease
the concentration of another substrate, depending on the concentrations of other proteins.
\begin{table}[htb]
  \makebox[\textwidth][c]{
  \small
    \begin{tabular}{cccccccccccc}
      Time & Cln3 & MBF & SBF &Cln1,2 & Cdh1 & Swi5 & Cdc14,20 & Clb5,6 & Sic1 & Clb1,2 & Mcm1/SFF\\
      \hline
      1&1&0&0&0&1&0&0&0&1&0&0\\
      2&0&1&1&0&1&0&0&0&1&0&0\\
      3&0&1&1&1&1&0&0&0&1&0&0\\
      4&0&1&1&1&0&0&0&0&0&0&0\\
      5&0&1&1&1&0&0&0&1&0&0&0\\
      6&0&1&1&1&0&0&0&1&0&1&1\\
      7&0&0&0&1&0&0&1&1&0&1&1\\
      8&0&0&0&0&0&1&1&0&0&1&1\\
      9&0&0&0&0&0&1&1&0&1&1&1\\
      10&0&0&0&0&0&1&1&0&1&0&1\\
      11&0&0&0&0&1&1&1&0&1&0&0\\
      12&0&0&0&0&1&1&0&0&1&0&0\\
      13&0&0&0&0&1&0&0&0&1&0&0\\
      \hline
    \end{tabular}
  }
  \caption{The temporal evolution of the Boolean cell-cycle model in \cite{Tang:2004}; corresponding to the biological cell-cycle sequence.}
  \label{table:timecourse}
\end{table}

Li \emph{et al.} \cite{Tang:2004} use their model to generate a time course of temporal
evolution of the cell-cycle network, shown in Table \ref{table:timecourse}.
This time course is in agreement with the behavior of the cell-division
process cycling through the four distinct phases $G_1$, $S$ (Synthesis), $G_2$, and $M$ (Mitosis).

We used this time course along with the network
as the only input to the algorithm to obtain all nested canalyzing
functions that interpolate the time course. In the fifth column of Table \ref{table:numFun}
we list the number of NCFs for each protein.
By requiring the Boolean function to be nested canalyzing, we have significantly reduced the number of possible functions
for each protein as it is evident when comparing the numbers in the third and fifth columns of Table \ref{table:numFun}.
However, even after this reduction, there are $330,559,488$ possible nested canalyzing models that fit the time course in Table \ref{table:timecourse}. To reduce the this number more, one needs to use additional time courses or request that the models
incorporate additional biological information about yeast cell cycle.

\begin{table}[htb]
\begin{center}
    \begin{tabular} {|l|c|c|c|c|}
      \hline
      Protein ($i$)& inputs & $f_i+I$ &  NCFs & NCFs in $f_i+I$\\
      \hline
      Cln3      & 1 & 1 & 2 & 1\\
      MBF       & 3 & 8 & 64 & 2\\
      SBF       & 3 & 8 & 64 & 2\\
      Cln1,2    & 1 & 1 & 2 & 1\\
      Cdh1      & 4 & 2048 & 736 & 12\\
      Swi5      & 4 & 2048 & 736 & 14\\
      Cdc20\&Cdc14& 3 & 8 & 64 & 4\\
      Clb5,6    & 3 & 8 & 64 & 3\\
      Sic1      & 5 & $2^{24}$ & 10,634 & 336\\
      Clb1,2    & 5 & $2^{24}$ & 10,634 & 61\\
      Mcm1/SFF  & 3 & 8 & 64 & 2\\
      \hline
    \end{tabular}
  \end{center}
  \caption{For each protein $i$, we list the number of inputs, the number of possible Boolean functions (the cardinality of $f_i+I$), the number of nested canalyzing functions with the given number of inputs, and finally the number of nested canalyzing functions in the model space $f_i+I$.}
  \label{table:numFun}
\end{table}
\subsection{Dynamics}
To analyze the dynamics of the resulting nested canalyzing models, we randomly sampled
2000 models and analyzed them. The average number of basins of attraction
(components) per network is
3.09, and the average size of the component containing the given trajectory is
1889. In only 6 models, the trajectory in Table \ref{table:timecourse} is not in the largest component,
however, the average size of the component containing the
trajectory is 833.5.

These results clearly show that nested canalyzing models for the cell cycle network
are in agreement with the original threshold model of Li \emph{et al.}, and since such models
are known to be robust and stable, any of these models
could be used as a model for the cell cycle in budding yeast. Furthermore, especially
when there is no evidence for choosing a particular type of functions, a
nested canalyzing function has an advantage over other possible choices.


\subsection{Comparison with Random Networks}
\begin{figure}[htb]
  \begin{center}
    \includegraphics[width=.5\textwidth]{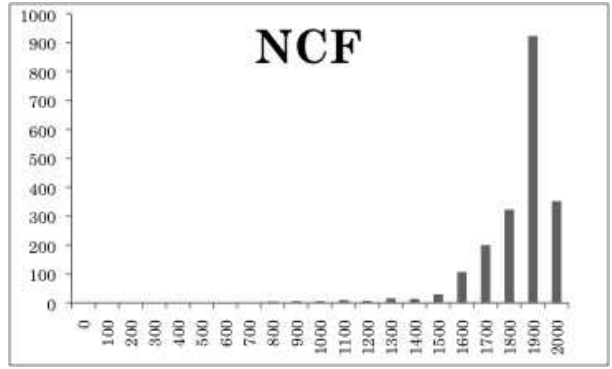}
    \caption{Nested canalyzing functions with the wiring diagram in Figure \ref{fig:dep} interpolating the time course in Table \ref{table:timecourse}.
      $x$-axis: size of basin of attraction for given trajectory,
      $y$-axis: number of networks observed, out of 2000.}
    \label{fig:ncf}
  \end{center}
\end{figure}

\begin{figure}[htb]
  \begin{center}
    \includegraphics[width=.5\textwidth]{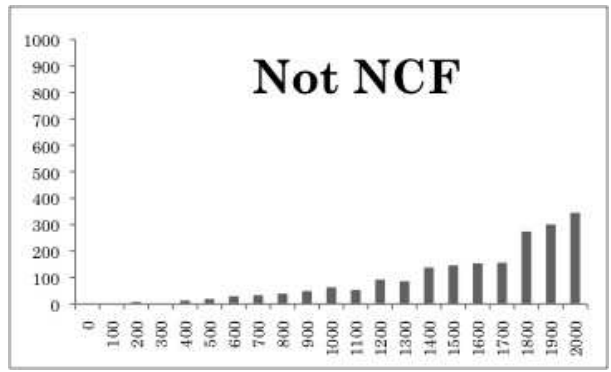}
    \caption{Not nested canalyzing functions with the wiring diagram in Figure \ref{fig:dep} interpolating the time course in Table \ref{table:timecourse}.
      $x$-axis: size of basin of attraction for given trajectory,
      $y$-axis: number of networks observed, out of 2000.}
    \label{fig:random}
  \end{center}
\end{figure}

To understand the effect of the network itself on its dynamics,
we sampled 2000 models on the same network where the local function of each
gene in each one of these models is chosen randomly from all possible functions in the model space.
 We found that the given
cell cycle trajectory has oftentimes much smaller basin of attraction,
and hence random functions on the cell cycle network could not in general produce
the desired dynamics.
A comparison of
the statistics from the sampled networks is shown in Figures \ref{fig:ncf} and
\ref{fig:random}.

\section{Conclusion}
In this paper we have presented an algorithm for identifying all Boolean nested canalyzing models
that fit a given time course or other input-output data sets.
Our algorithm uses methods from computational algebra to present the model
space as an algebraic variety. The intersection of this variety with the
variety of all NCFs, which was parameterized in \cite{AB07}, gives us the set of all
NCFs that fit the data. We demonstrate our algorithm by finding all nested canalyzing
models of the cell cycle network form Li \emph{et al.} \cite{Tang:2004}.
We then showed that the dynamics of almost any of these models is strikingly similar
to that of the original threshold model.
Unless the chosen model is required to meet other conditions, and in that case the model
space will be reduced further, any one of the models that our algorithm found is an acceptable
model of the cell cycle process in the Budding yeast.

One limitation of the current algorithm,
which we left for future work,
is that it does not distinguish between activation and inhibition in the network
as we do not have a systematic method of knowing when a given variable in a (nested canalyzing) polynomial is an activator or inhibitor.

As our algorithm relies heavily on different Gr\"obner based computations,
the current implementation in Singular allows a given gene to have at most 5 regulators.
This is due to the fact that the number of monomials then is 32 which is already a burden
especially when the primary decomposition of an ideal is what we are after.
We are working on a better implementation so that we can infer larger and denser networks.


\bibliographystyle{plain}
\bibliography{./submitted}

\begin{thebibliography}{10}

\bibitem{AdamsLost}
W.~Adams and P.~Loustaunau.
\newblock {\em An introduction to {G}r\"obner bases}, volume~3 of {\em Graduate
  Studies in Mathematics}.
\newblock American Mathematical Society, Providence, RI, 1994.

\bibitem{singular}
W.~Decker, G.-M. Greuel, G.~Pfister, and H.~Sch{\"o}nemann.
\newblock {\sc Singular} {3-1-1} --- {A} computer algebra system for polynomial
  computations, 2010.
\newblock http://www.singular.uni-kl.de.

\bibitem{DJLS}
E.~Dimitrova, A.~Jarrah, R.~Laubenbacher, and B.~Stigler.
\newblock A {G}r\"obner fan method for biochemical network modeling.
\newblock In {\em I{SSAC} 2007}, pages 122--126. ACM, New York, 2007.

\bibitem{DVML}
E.~Dimitrova, J.~McGee, R.~Laubenbacher, and P.~Vera~Licona.
\newblock Comparison of discretization methods for network inference.
\newblock {\em Journal of Computational Biology}, 2010.
\newblock In Press.

\bibitem{chain}
I.~Gat-Viks and R.~Shamir.
\newblock Chain functions and scoring functions in genetic networks.
\newblock {\em Bioinformatics}, 19:108--117, 2003.

\bibitem{harris}
S.~Harris, B.~Sawhill, A.~Wuensche, and S.~Kauffman.
\newblock A model of transcriptional regulatory networks based on biases in the
  observed regulation rules.
\newblock {\em Complex.}, 7(4):23--40, 2002.

\bibitem{Hinkelmann}
F.~Hinkelmann.
\newblock Singular implementation of {NCF} inferring algorithm.
\newblock Available at http://www.math.vt.edu/people/fhinkel/ncf.lib, 2010.

\bibitem{AB07}
A.~Jarrah and R.~Laubenbacher.
\newblock Discrete models of biochemical networks: The toric variety of nested
  canalyzing functions.
\newblock In H.~Anai, K.~Horimoto, and T.~Kutsia, editors, {\em Algebraic
  Biology}, number 4545 in LNCS, pages 15--22. Springer, 2007.

\bibitem{JLSS}
A.~Jarrah, R.~Laubenbacher, B.~Stigler, and M.~Stillman.
\newblock Reverse-engineering of polynomial dynamical systems.
\newblock {\em Advances in Applied Mathematics}, 39:477--489, 2007.

\bibitem{JLR}
A.~Jarrah, B.~Raposa, and R.~Laubenbacher.
\newblock Nested canalyzing, unate cascade, and polynomial functions.
\newblock {\em Physica D}, 233:167--174, 2007.

\bibitem{Kauffman:03}
S.~Kauffman, C.~Peterson, B.~Samuelsson, and C.~Troein.
\newblock Random boolean network models and the yeast transcriptional network.
\newblock {\em PNAS}, 100(25):14796--14799, 2003.

\bibitem{Kauffman:04}
S.~Kauffman, C.~Peterson, B.~Samuelsson, and C.~Troein.
\newblock {Genetic networks with canalyzing Boolean rules are always stable}.
\newblock {\em PNAS}, 101(49):17102--17107, 2004.

\bibitem{Kauffman:69}
S.~A. Kauffman.
\newblock Metabolic stability and epigenesis in randomly constructed genetic
  nets.
\newblock {\em Journal of Theoretical Biology}, 22:437--467, 1969.

\bibitem{LS}
R.~Laubenbacher and B.~Stigler.
\newblock A computational algebra approach to the reverse-engineering of gene
  regulatory networks.
\newblock {\em Journal of Theoretical Biology}, 229:523--537, 2004.

\bibitem{Tang:2004}
F.~Li, T.~Long, Y.~Lu, Q.~Ouyang, and C.~Tang.
\newblock The yeast cell-cycle network is robustly designed.
\newblock {\em PNAS}, 101:4781, 2004.

\bibitem{david}
D.~Murrugarra.
\newblock Multi-states nested canlyzing functions.
\newblock 2010.
\newblock preprint.

\bibitem{wilhelm_2007}
S.~Nikolajewaa, M.~Friedela, and T.~Wilhelm.
\newblock Boolean networks with biologically relevant rules show ordered
  behavior.
\newblock {\em Biosystems}, 90(1):40--47, 2007.

\bibitem{ode-unate}
R.~Porreca, E.~Cinquemani, J.~Lygeros, and G.~Ferrari-Trecate.
\newblock {Identification of genetic network dynamics with unate structure}.
\newblock {\em Bioinformatics}, 26(9):1239--1245, 2010.

\bibitem{Raey_2002}
L.~Raeymaekers.
\newblock Dynamics of boolean networks controlled by biologically meaningful
  functions.
\newblock {\em Journal of Theoretical Biology}, 218(3):331 -- 341, 2002.

\bibitem{shmul_2003}
I.~Shmulevich, H.~L\"{a}hdesm\"{a}ki, E.~R. Dougherty, J.~Astola, and W.~Zhang.
\newblock {The role of certain Post classes in Boolean network models of
  genetic networks}.
\newblock {\em PNAS}, 100(19):10734--10739, 2003.

\bibitem{reveng-review}
C.~Sima, J.~Hua, and S.~Jung.
\newblock Inference of gene regulatory networks using time-series data: A
  survey.
\newblock {\em Current Genomics}, 10(14):416--429, 2009.

\bibitem{Steggles}
L.~J. Steggles, R.~Banks, O.~Shaw, and A.~Wipat.
\newblock Qualitatively modelling and analysing genetic regulatory networks: a
  {P}etri net approach.
\newblock {\em Bioinformatics}, 23:336--343, 2007.

\bibitem{SJLS}
B.~Stigler, A.~Jarrah, M.~Stillman, and R.~Laubenbacher.
\newblock Reverse-engineering of dynamic networks.
\newblock {\em Annals of the New York Academy of Sciences}, 1115:168--177,
  2007.

\bibitem{LM_thomas}
R.~Thomas.
\newblock Boolean formalisation of genetic control circuits.
\newblock {\em Journal of Theoretical Biology}, 42:565--583, 1973.

\bibitem{V-CJL}
A.~Veliz-Cuba, A.~Jarrah, and R.~Laubenbacher.
\newblock Polynomial algebra of discrete models in systems biology.
\newblock {\em Bioinformatics}, 26(13):1637--1643, 2010.

\bibitem{wad}
C.~H. Waddington.
\newblock Canalisation of development and the inheritance of acquired
  characters.
\newblock {\em Nature}, 150:563--564, 1942.

\end{thebibliography}
\end{document}